\documentclass[11pt]{article}

\usepackage{latexsym}
\usepackage{amsfonts}
\usepackage{amsfonts}
\usepackage{enumerate}
\usepackage{multicol}
\usepackage{graphicx}
\usepackage{amssymb}
\usepackage{amsmath}
\usepackage{amsthm}
\usepackage{epic}
\usepackage{sidecap}
\usepackage{tikz}

\begin{document}
\parindent=0.2in
\parskip 0in

\title{Report on the 12th Annual USA Junior Mathematical Olympiad \\[.4in]}

\author{B\'{e}la Bajnok and Evan Chen} 

\date{September 10, 2021}

\maketitle

\vspace*{.2in}

The USA Junior Mathematical Olympiad (USAJMO) is the final round in the American Mathematics Competitions series for high school students in grades 10 or below, organized each year by the Mathematical Association of America.  The competition follows the style of the International Mathematics Olympiad: it consists of three problems each on two consecutive days, with an allowed time of four and a half hours both days.

The 12th annual USAJMO was given on Tuesday, April 13, 2021 and Wednesday, April 14, 2021.  This year, 259 students were invited to take the USAJMO and, as in 2020, the competition was administered online.  The names of winners and honorable mentions, as well as more information on the American Mathematics Competitions program, can be found on the site https://www.maa.org/math-competitions.  Below we present the problems and solutions of the competition; a similar article for the USA Mathematical Olympiad (USAMO) can be found in Mathematics Magazine.

The problems of the USAJMO are chosen -- from a large collection of proposals submitted for this purpose -- by the USAMO/USAJMO Editorial Board, whose co-editors-in-chief this year were Evan Chen and Jennifer Iglesias, with associate editors Ankan Bhattacharya, John Berman, Zuming Feng, Sherry Gong, Alison Miller, Maria Monks Gillespie, and Alex Zhai.  This year's problems were created by Ankan Bhattacharya, Vincent Huang, Mohsen Jamaali, Carl Schildkraut, Brandon Wang, and Alex Zhai.

The solutions presented here are those of the present authors, relying in part on the submissions of the problem authors and members of the editorial board.
Each problem was worth 7 points; the nine-tuple $(n; a_7, a_6, a_5, a_4, a_3, a_2, a_1, a_0)$ states the number of students who submitted a paper for the relevant problem, followed by the numbers who scored $7, 6, \dots, 0$ points, respectively.

\vspace{.1in}

\noindent {\bf Problem 1} $(236;83,23,4,0,5,9,26,86)$; {\em proposed by Vincent Huang}.   Let $\mathbb N$ denote the set of positive integers.
	Find all functions $f \colon \mathbb N\to\mathbb N$
	such that for all positive integers $a$ and $b$,
	\[ f(a^2+b^2)=f(a)f(b) \quad\text{and}\quad f(a^2)=f(a)^2. \]

\vspace{.1in}

\noindent {\em Solution}.  It is easy to see that the constant function $f \equiv 1$ satisfies the two given properties; we will use induction to prove that there are no others.  Since, according to the second condition, $f(1^2)=f(1)^2$, we see that $f(1)=1$; and now from the first condition we get $f(2)=f(1^2+1^2)=f(1)^2=1$ as well.  For the inductive step, we let $n \geq 3$ be any integer, and suppose that for every positive integer $m < n$ we have $f(m)=1$; we need to prove that $f(n)=1$ holds also.

We start with the identity (familiar from Pythagorean triples) that for all positive integers $u$ and $v$ with $u > v$ we have
$$(u^2-v^2)^2 + (2uv)^2 = (u^2+v^2)^2,$$
and therefore, using the two given conditions, we find that
$$f(u^2-v^2)f(2uv)=f \left((u^2+v^2)^2 \right) = \left(f(u^2+v^2) \right)^2 = \left(f(u)f(v)\right)^2.$$
In particular, if $u>v$ are positive integers for which $f(u)=f(v)=1$, then $f(u^2-v^2)=1$ and $f(2uv)=1$.

We consider two cases: when $n$ is odd and when $n$ is even.  If $n=2k+1$ for some positive integer $k$, then, since $f(k+1)=f(k)=1$ by our inductive assumption, we get $f((k+1)^2-k^2)=f(n)=1$.

Similarly, when $n=2k$, we have $f(k)=f(1)=1$, and thus $f(2k)=f(n)=1$.  Our proof is now complete.

\vspace{.1in}

\noindent {\bf Problem 2} $(176;66,1,1,0,0,0,8,100)$; {\em proposed by Ankan Bhattacharya}.     Rectangles $BCC_1B_2$, $CAA_1C_2$, and $ABB_1A_2$ are erected outside an acute triangle $ABC$. Suppose that
	\[ \angle BC_1C + \angle CA_1A + \angle AB_1B = 180^\circ. \]
Prove that lines $B_1C_2$, $C_1A_2$, and $A_1B_2$ are concurrent.

\vspace{.1in}

\noindent {\em Solution}.  Let $\omega_A$, $\omega_B$, and $\omega_C$ be the circumcircles of rectangles $BCC_1B_2$, $CAA_1C_2$, and $ABB_1A_2$, respectively.

\begin{figure}[ht]
\centering
\includegraphics[scale=.8]{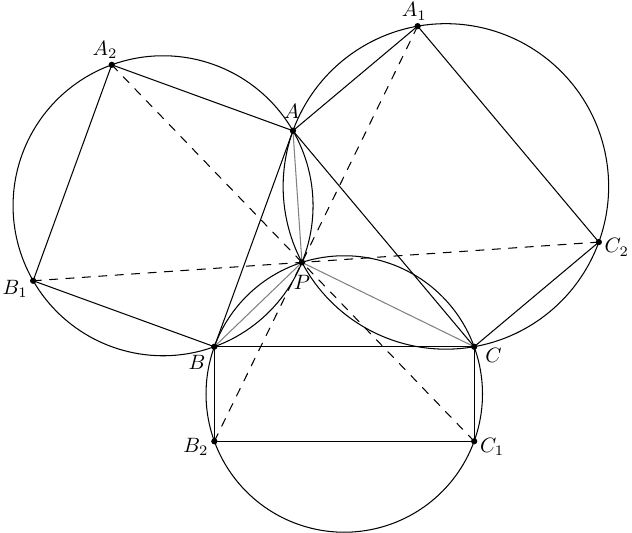}
\caption{The main circles used in Problem 2}
\end{figure}
  Define $P$ to be the foot of the altitude from $A$ to $B_1C_2$.

Observe that, since $\angle ACC_2$ and $\angle APC_2$ are both $90^\circ$, $P$ lies on $\omega_B$ by the inscribed angle theorem and, therefore, $\angle APC$ and $\angle CA_1A $ are supplementary angles.  A similar argument shows that $P$ lies on $\omega_C$ and thus $ \angle APB$ and $\angle AB_1B$ are supplementary angles as well.

But then
\[ \angle BPC = 360^\circ - (\angle APC + \angle APB) =  \angle CA_1A + \angle AB_1B,\]
which, by the given equation, yields
$$ \angle BPC = 180^\circ - \angle BC_1C, $$ and thus $P$ lies on $\omega_A$ as well.  Therefore, $P$ is the (unique) common point of all three circles.

Similar arguments would prove that the feet of the altitudes from $B$ and $C$ to $C_1A_2$ and $A_1B_2$, respectively, are on each of the three circles, and thus must coincide with $P$.  But then lines $B_1C_2$, $C_1A_2$, and $A_1B_2$ are concurrent, as claimed.

\vspace{.1in}

\begin{figure}[ht]
	\centering
	\includegraphics[scale=0.8]{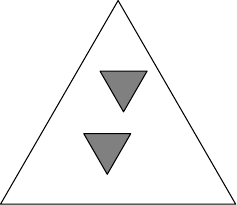}
	\caption{An example for Problem 3 with two triangles}
	\label{JMO 3 statement}
\end{figure}

\noindent {\bf Problem 3} $(159;13,11,2,0,1,6,26,100)$; {\em proposed by Alex Zhai}.    An equilateral triangle $\Delta$ of side length $L > 0$ is given.
	Suppose that $n$ equilateral triangles with side length $1$
	and with non-overlapping interiors are drawn inside $\Delta$,
	such that each unit equilateral triangle has sides parallel to $\Delta$,
	but with opposite orientation.
	(An example with $n=2$ is drawn in Figure \ref{JMO 3 statement}.)

Prove that $n \le \frac{2}{3} L^2. $

\vspace{.1in}

\begin{figure}[ht]
	\centering
	\includegraphics[scale=0.8]{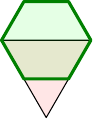}
	\caption{The key hexagon in the solution of Problem 3}
	\label{hexagon fig}
\end{figure}

\noindent {\em Solution}.  For each triangle, we draw a green regular hexagon
of side length $1/2$ as shown in Figure \ref{hexagon fig}.  Because the edges of $\Delta$ slope upwards at the same angle
as the edges in the upper half of the green hexagons,
it follows that the green hexagons lie entirely inside $\Delta$.

Our main claim is that the green hexagons are pairwise disjoint.  To see this,
anchor each triangle by the midpoint of the top side.
By fixing one anchor and moving another anchor as close as possible
so that the triangles do not overlap,
we find that the boundary of the allowable region
for the second anchor is a hexagon of side length $1$ centered
at the first anchor.
This is equivalent to the hexagons not overlapping.

Now since each hexagon has area $\frac{3\sqrt{3}}{8}$ and
lies inside $\Delta$, and the total area of $\Delta$ is $\frac{\sqrt3}{4}L^2$,
we conclude that
\[ \frac{3\sqrt3}{8} \cdot n \le \frac{\sqrt3}{4} L^2,\]
and thus $n \le \frac 23 L^2,$ as claimed.

\begin{figure}[ht]
	\centering
	\includegraphics{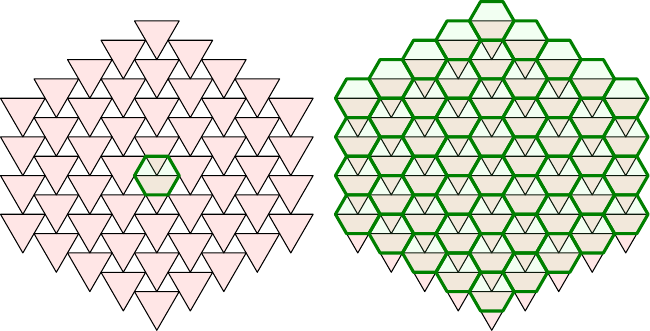}
	\caption{Tessellation achieving almost $\frac23$ density.
	(In the illustration on the left,
	only one of the green hexagons is drawn in.
	The version on the right has all the green hexagons.)}
	\label{fig:tessellation}
\end{figure}

{\em Remark}.
The constant $\frac23$ cannot be improved.
The tessellation in Figure~\ref{fig:tessellation} illustrates how one can achieve $n \ge \left( \frac 23 - \epsilon \right) L^2$,
for any $\epsilon > 0$.

\vspace{.1in}

\noindent {\bf Problem 4} $(228;56,28,6,5,9,7,34,83)$; {\em proposed by Brandon Wang}.
	Carina has three pins, labeled $A$, $B$, and $C$, respectively,
	located at the origin of the coordinate plane.
	In a \emph{move}, Carina may move a pin to
	an adjacent lattice point at distance $1$ away.
	What is the least number of moves that Carina can make
	in order for triangle $ABC$ to have area $2021$?

	(A lattice point is a point $(x,y)$ in the coordinate plane
	where $x$ and $y$ are both integers, not necessarily positive.)

\vspace{.1in}

\noindent {\em Solution}. The answer is $128$.
One construction with $128$ moves is given by moving the pins to the points
\begin{align*}
	A &= (-3,-18), \\
	B &= (61,0), \\
	C &= (0,46).
\end{align*}
We can verify that the number of moves is
\[3+18+61+0+0+46=128,\]
and that the area of the triangle (using, for example, the shoelace formula) is
\[ \frac{61 \cdot 46 + 18 \cdot 61 + 3 \cdot 46}{2} = 2021. \]

\begin{figure}[ht]
	\centering
	\includegraphics{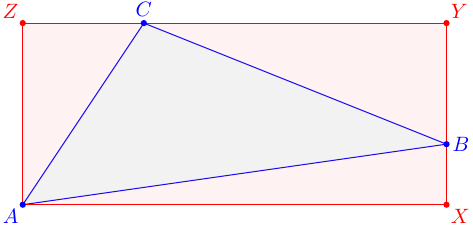}
	\caption{The bounding box of a triangle}
	\label{fig J4}
\end{figure}

We now prove that at least $128$ moves are required.
Define the {\em bounding box} of triangle $ABC$
to be the smallest axis-parallel rectangle which
contains all three of the vertices $A$, $B$, $C$;
an example is illustrated in Figure~\ref{fig J4}.

We first establish the (well known) claim that the area of a triangle is at most half the area of its bounding box.
Note that all four sides of the bounding box must contain at least one vertex of the triangle, which implies that at least one vertex, say $A$, is at a corner.  We label the other three corners of the rectangle by $X$, $Y$, and $Z$, so that $Y$ is opposite $A$.  If $B=Y$ or $C=Y$, then $\triangle ABC$ is contained in $\triangle AXY$ or $\triangle AZY$, so our claim obviously holds.  
If $B \neq Y$ and $C \neq Y$, then these two vertices lie on different sides of the line $AY$, and thus one of them lies on $XY$ and the other on $YZ$; we assume that $B$ is on $XY$ and $C$ is on $YZ$.  Note that $Z$ is at least as far from the line $AB$ as $C$ is, which means that the area of $\triangle ABC$ is at most as much as the area of $\triangle ABZ$, implying our claim.     

We can now complete our proof, as follows.  Suppose that $n$ moves are made,
leading to a bounding box of width $w$ and length $\ell$.  Since any move changes the dimensions
of the bounding box by at most $1$, we have $\ell + w \leq n$, and thus for the area of triangle $ABC$ we have
$$2021 \leq \frac12 \ell w = \frac18 \left( (\ell +w)^2 -  (\ell -w)^2 \right) \leq \frac18  (\ell +w)^2 \leq \frac18 n^2,$$
which yields $n \ge 128$, as claimed.

\vspace{.1in}

\noindent {\bf Problem 5} $(202;42,11,8,2,2,31,13,93)$; {\em proposed by Carl Schildkraut}.   A finite set $S$ of positive integers has the property that,
	for each $s\in S$, and each positive integer divisor $d$ of $s$,
	there exists a unique element $t\in S$ satisfying $\gcd(s,t) = d$.
	(The elements $s$ and $t$ could be equal.)

	Given this information, find all possible values for the
	number of elements of $S$.

\vspace{.1in}

\noindent {\em Solution.}  We claim that the possible sizes are 0 and the nonnegative integer powers of 2.  Since $S=\emptyset$ and $S=\{1\}$ obviously work, we need to show that a set $S$ of size $n \geq 2$ satisfying the requirements exists if, and only if, $n=2^k$ for some positive integer $k$.

We start by verifying that these values are indeed possible.  For a given positive integer $k$, we construct a set of size $2^k$ as follows.  Suppose that $p_1, q_1, p_2, q_2, \dots, p_k, q_k$ are $2k$ pairwise distinct positive primes; for an ordered pair of subsets $(I,J)$ of $[k]=\{1,2, \dots, k\}$, we will use the notation
$$s(I,J)=\prod_{i \in I} p_i \cdot \prod_{j \in J} q_j.$$  (Recall that the empty product equals 1.)
We then consider
the set $$S=\{s(I,J) \mid I \subseteq [k], J=[k] \setminus I\}.$$  Since the elements of $S$ are then in a bijection with the subsets $I$ of $[k]$, we see that $|S|=2^k$.  We need to show that $S$ satisfies the required property.

Given an element $s(I,J)$ of $S$, we see that its positive divisors are of the form $s(I_0,J_0)$ where $I_0 \subseteq I$ and $J_0 \subseteq J$.  Note also that, since $I$ and $J$ form a partition of $[k]$, the sets $I'=I_0 \cup (J \setminus J_0)$ and $J'=J_0 \cup (I \setminus I_0)$ do as well, and thus
$s(I',J')$ is an element of $S$; in fact it is the unique element of $S$ whose greatest common divisor with $s(I,J)$ equals $s(I_0,J_0)$.   Therefore, the set $S$ we constructed satisfies the requirement of the problem.

It remains to be shown that if $S$ is a set satisfying the property and it has size $n  \geq 2$, then $n=2^k$ for some positive integer $k$.  Let $s \geq 2$ be any element of $S$, and let $p$ be any positive prime divisor of $s$.  We can then write $s=p^e \cdot u$ for some positive integers $e$ and $u$ where $u$ is not divisible by $p$.  We claim that $e=1$.

Denoting by $d(m)$ the number of positive divisors of a positive integer $m$, we have $d(s)=(e+1) \cdot d(u)$; in fact, $s$ has exactly $e  \cdot d(u)$ positive divisors that are divisible by $p$ and $d(u)$ that are not.  By our assumption, there is an element $t$ of $S$ for which $\gcd (t,s)=p$.  Let us assume that $e \geq 2$; we can then see that $t=p \cdot v$ for some positive integer $v$ that is not divisible by $p$.  Furthermore, $d(t)=2d(v)$, and $t$ has exactly $d(v)$ positive divisors that are divisible by $p$ and also $d(v)$ that are not.

Now according to our requirement for $S$, the positive divisors of any element $s$ of $S$ are in a one-to-one correspondence with the elements of $S$ via the map $a \mapsto \gcd(a,s)$, and thus $|S|=d(s)$.  (In particular, all elements of $S$ must have the same number of positive divisors.)  Furthermore, for any prime divisor $p$ of $s$, we have $p|a$ if and only if $p| \gcd(a,s)$.  Therefore, $S$ has exactly $e  \cdot d(u)$ elements that are divisible by $p$ and $d(u)$ that are not.  With the same reasoning, $S$ has exactly $d(v)$ elements that are divisible by $p$ and $d(v)$ that are not.  But then $d(u)=d(v)$ and $e  \cdot  d(u)=d(v)$, from which $e=1$.

This establishes the fact that each element of $S$ is a product of the same number of pairwise distinct prime numbers.  If $s \in S$ is the product of $k$ distinct primes, then $|S|=d(s)=2^k$.  This completes our proof.

\vspace{.1in}

\noindent {\bf Problem 6} $(147;16,2,0,0,4,21,2,102)$; {\em proposed by Mohsen Jamaali}.   Let $n \ge 4$ be an integer.
	Find all positive real solutions to the following
	system of $2n$ equations:
	\begin{align*}
		a_1 = &\frac{1}{a_{2n}} + \frac{1}{a_{2}}, & a_2 &= a_1 + a_3, \\[1ex]
		a_3 = &\frac{1}{a_{2}} + \frac{1}{a_{4}}, & a_4 &= a_3 + a_5, \\[1ex]
		a_5 = &\frac{1}{a_{4}} + \frac{1}{a_{6}}, & a_6 &= a_5 + a_7, \\[1ex]
& \vdots & & \vdots \\
		a_{2n-1} = &\frac{1}{a_{2n-2}} + \frac{1}{a_{2n}}, & a_{2n} &= a_{2n-1} + a_1.
	\end{align*}

\vspace{.1in}

\noindent {\em First solution.}  It is easy to see that a solution is provided by $a_1=a_3=\cdots=a_{2n-1}=1$ and $a_2=a_4=\cdots=a_{2n}=2$; we prove that there are no others.

Taking indices modulo $2n$ and eliminating terms of odd indices, for each $i=1,2,\dots,n$ we have
\begin{eqnarray}
a_{2i}=\frac{1}{a_{2i-2}}+\frac{2}{a_{2i}}+\frac{1}{a_{2i+2}}.
\end{eqnarray}
Adding up these equations yields
\begin{eqnarray}
\sum_{i=1}^n a_{2i}=\sum_{i=1}^n \frac{4}{a_{2i}}.
\end{eqnarray}
According to the harmonic mean--arithmetic mean inequality,
\begin{eqnarray}
\frac{n}{\sum_{i=1}^n \frac{1}{a_{2i}}} \leq \frac{\sum_{i=1}^n a_{2i}}{n},
\end{eqnarray}
so by (2) we get $\sum_{i=1}^n a_{2i} \geq 2n$.

Now dividing both sides of (1) by $a_{2i}$ yields
\begin{eqnarray}
1=\frac{1}{a_{2i-2}a_{2i}}+\frac{2}{a_{2i}^2}+\frac{1}{a_{2i}a_{2i+2}},
\end{eqnarray}
and adding the equations results in
\begin{eqnarray}
n=\sum_{i=1}^n \left( \frac{1}{a_{2i}}+ \frac{1}{a_{2i+2}} \right)^2.
\end{eqnarray}
Now we use the quadratic mean--arithmetic mean inequality, which gives
\begin{eqnarray}
\frac{1}{n} \cdot \sum_{i=1}^n \left( \frac{1}{a_{2i}}+ \frac{1}{a_{2i+2}} \right)^2 \geq \left( \frac{\sum_{i=1}^n \left( \frac{1}{a_{2i}}+ \frac{1}{a_{2i+2}} \right)}{n} \right)^2,
\end{eqnarray}
so by (5) and (2) we get
$$1 \geq \left( \frac{ \sum_{i=1}^n \frac{2}{a_{2i}}}{n} \right)^2 =  \left( \frac{ \frac{1}{2} \cdot  \sum_{i=1}^n a_{2i} }{n} \right)^2 $$
and thus $\sum_{i=1}^n a_{2i} \leq 2n$.

This means that each of our inequalities is an equality, and therefore $a_{2i}=2$ for all $i$.  This in turn implies that $a_{2i-1}=1$ for all $i$, as claimed.

\vspace{.1in}

\noindent {\em Second solution.}  We write $m= \min \{a_{2i} \mid 1 \leq i \leq n\}$ and $M= \max \{a_{2i} \mid 1 \leq i \leq n\}$, and assume that $m=a_{2j}$ and $M=a_{2k}$.  Then equation (1) from the first solution yields
$$m=\frac{1}{a_{2j-2}}+\frac{2}{m}+\frac{1}{a_{2j+2}} \geq \frac{1}{M}+\frac{2}{m}+\frac{1}{M}$$ and
$$M=\frac{1}{a_{2k-2}}+\frac{2}{M}+\frac{1}{a_{2k+2}} \leq \frac{1}{m}+\frac{2}{M}+\frac{1}{m}.$$
Therefore, $$m \geq \frac{2}{m}+\frac{2}{M} \geq M,$$ which can only occur when $m=M$.  Therefore, all $a_{2i}$ are equal, from which  $a_2=a_4=\cdots=a_{2n}=2$ and $a_1=a_3=\cdots=a_{2n-1}=1$, as claimed.

\vspace*{.4in}

{\bf B\'ela Bajnok} (AMCDirector@MAA.edu, MR author ID 314851, ORCID number 0000-0002-9498-1596) is the director of the American Mathematics Competitions program of the MAA.  He is also a professor of mathematics at Gettysburg College, and the author of two books: {\em An Invitation to Abstract Mathematics}, now in second edition by Springer, and {\em Additive Combinatorics: A Menu of Research Problems}, published by CRC Press.

\vspace*{.2in}

	{\bf Evan Chen} (\texttt{evan@evanchen.cc}, MR author ID 1158569, ORCID number 0000-0001-9550-5068) is a mathematics Ph.D. student at the Massachusetts Institute of Technology.  He is the Co-Editor-in-Chief of the USA(J)MO Editorial Board and a coach of the USA team at the International Mathematical Olympiad.  His book {\em Euclidean Geometry in Mathematical Olympiads} was published by the MAA.

\end{document}